\documentclass[10pt]{article}
\usepackage[all]{xy}
\usepackage{lscape}
\usepackage[ansinew]{inputenc}
\usepackage[psamsfonts]{amssymb}
\def\bk{I\!\!k} 

\newcommand{\sss}[0]{Serre spectral sequence}

\title{ The loop-coproduct spectral sequences.}
\author{ Jean-Fran\c cois LE BORGNE }

\begin{document}

\maketitle

\begin{abstract} 
Let $M$ be a closed oriented $d$-dimensionnal manifold and let $LM$
be the
space of free loops on $M$. 
In this paper, we give a geometrical interpretation of the loop-coproduct
and we study it's compatibility with the Serre spectral
sequence associated to the fibration  $\Omega M \to LM
\stackrel{ev(0)}{\to} M$.  Then, we show that the spectral
sequence associated to the free loop fibration $LN \to LX \to LM$ of
some Serre fibration $N \to X \to M$ is a spectral sequence of Frobenius algebra.

\end{abstract}

\vspace{12mm}

\vspace{5mm}\noindent {\bf AMS Classification} : 55P35, 54N45,55N33, 17A65,

81T30, 17B55

\vspace{12mm}\noindent {\bf Key words} : free loop space, loop-homology,
Serre spectral sequence, loop-coproduct.

\vspace{5mm}

\centerline{\bf \Large  Introduction.}

\vspace{5mm}

Let $M$ be a $1$-connected closed oriented $d$-manifold and let $LM$
be the Hilbert manifold homotopically equivalent to $M^{S^1}$ the
space of free loops of $M$ \cite{CHA3}. We denote by $ev(0)$
(respectively ev(1/2)) the evaluation map at $0$ (respectively $1/2$).
$$ev(0):LM \to M \qquad \gamma \mapsto \gamma(0)$$
$$ev(1/2):LM \to M \qquad \gamma \mapsto \gamma(1/2)$$

We put $$L_{1/2}M=:\{\gamma \in LM  \quad /  \quad ev(1/2)(\gamma)=ev(0)(\gamma)
\}$$ and $$\Omega_{1/2} M =:\{\gamma \in \Omega M  \quad /  \quad ev(1/2)= * \}$$ where $*$
denotes the base point of $M$.
Then we have the pull back diagrams:
$$
\xymatrix{
L_{1/2}M \ar[r]^{\tilde{j}} \ar[d]^{ev(1/2)} & LM \ar[d]^{ev(0) \times
  ev(1/2)}\\
M \ar[r]^{\Delta} & M \times M }
\xymatrix{
\Omega_{1/2} M \ar[r]^{\tilde{i}} \ar[d]^{ev(1/2)} & \Omega M \ar[d]^{ev(1/2)}\\
\ast \ar[r]^{i} & M }
$$
where $\Delta$ is the diagonnal embedding and $i$ the canonical
inclusion of the base point of $M$. This diagrams define the
$d$-codimensionnal embeddings $\tilde{j}$ and $\tilde{i}$.
We need again the following notations:
$$LM \times_M LM = \{(\alpha , \beta) \in LM \times LM  \quad /  \quad \alpha(0)
= \beta(0) \}.$$ 
We have the following commutative pull-back diagram \cite{CJY3}:

$$
\xymatrix{
LM \times_M LM \ar[r]^{\tilde{\Delta}} \ar[d]^{ev_{\infty}} & LM \times
LM \ar[d]^{ev(0) \times ev(0)}\\
M \ar[r]^{\Delta} & M \times M}
$$
where $ev_{\infty}$ denotes the evaluation at $0$: $$LM \times_M LM
\to M \qquad (\alpha , \beta ) \mapsto \alpha(0) = \beta(0) $$.

We define
$$\gamma: LM \times_M LM \to L_{1/2}M $$
$$(\alpha ,  \beta) \mapsto \alpha * \beta$$
with $\alpha * \beta = \alpha(2t)$ if $t \in [0 , 1/2]$ and $\beta(2t-1)$
if $t \in [1/2,1]$ the reparametrisation of a couple of same basepoint
loops.
We remark that $\gamma$ is an homeomorphism and we have 
$$\gamma^{-1}: L_{1/2}M \to LM \times_M LM $$
$$\eta \mapsto (\eta(t/2),\eta(1/2+t/2))$$
We denote $\gamma_{\omega}$ (resp. $\gamma_{\omega}^{-1}$) the
restriction of $\gamma$ (resp. $\gamma^{-1}$) to $\Omega M \times
\Omega M $ (resp. $\Omega_{1/2} M$). 
We denote by $comp: LM \times_M LM \to LM $ the composition of free
loops and $comp_{\omega}$ its restriction to pointed loops. Remark
that $comp = \tilde{j} \circ \gamma$ for free loops and $comp_{\omega}=\tilde{i}
\circ \gamma_{\omega}$ for pointed loops.

Assume that $\bk$ is a fixed ring. The following homology groups are
assumed with coefficients in $\bk$.
M.Chas and D.Sullivan have constructed a product $$P:H_*(LM \times LM)
\to H_{*-d}(LM)$$ called the loop product such that the desuspended
homology of $LM$ namely $\mathbb{H}_*(LM) =: H_{*+d}(LM)$ is a
commutative graded algebra. With our notations, $P=:comp_* \circ
\tilde{\Delta}_!$.

In \cite{CG3}, R.Cohen and V.Godin have constructed a coproduct 
$$\Phi: H_*(LM) \to H_{*-d}(LM \times LM).$$ 
Since $\tilde{j} : L_{1/2}M
\hookrightarrow LM$ is a smooth finite codimensionnal embedding of
Hilbert manifolds, we can define $\tilde{j}_! : H_*(LM) \to
H_{*-d}(L_{1/2}M)$ (see \cite{CJY3}, \cite{JF3} or \cite{SAD3} for
details).
Furthermore, considering $\gamma$ as a $0$-codimensionnal smooth
embedding, we can also consider $$\gamma_!: H_*(L_{1/2}M) \to H_*(LM
\times_M LM).$$

\noindent {\bf Remark:} The coproduct on based loop space $comp_{{\omega}_!}$ is
the coproduct defined by Sullivan in \cite{SUL3}.

The remaining of this paper consist to to prove the following results:

\noindent {\bf Proposition 1:} 
{\it Let $h: X \to Y$ a diffeomorphism between
  smooth Hilbert closed connected manifolds without boundary. We have: 
$h_! = h^{-1}_*$}

\noindent {\bf Corollary:} 
{\it $\gamma_!=\gamma_*^{-1}$ and ${\gamma_{\omega}}_!={\gamma_{\omega}}_*^{-1}$}

So that we can define $comp_! =: \gamma_! \circ \tilde{j}_! = \gamma^{-1}_*
\circ \tilde{j}_!$ (by Proposition 1). As the same, we define
${comp_{\omega}}_! =: {\gamma_{\omega}}_! \circ \tilde{i}_! =
{\gamma_{\omega}^{-1}}_* \circ \tilde{i}_! $.
From the definition of $\Phi$, we obtain immediately:

\noindent {\bf Theorem 1:} 
{\it According to the preceding notations, $\Phi =
  \tilde{\Delta}_* \circ comp_!$ .}

Remark that the loop-coproduct $\Phi= \tilde{\Delta}_* \circ comp_!$ can be thincked
as Poincar\'e-dual of the loop-product $P = comp_* \circ
\tilde{\Delta}_! $.

Now, let us consider $\{E^*_{*,*}[ev(0)] \}$ the \sss $\quad$ associated to the fibration $\Omega M \to LM
\stackrel{ev(0)}{\to} M$.

\noindent {\bf Theorem 2:} 
{\it The spectral sequence  $\{E^*_{*,*}[ev(0)] \}$ is
  comultiplicative and converges to the coalgebra $(H_*(LM),
  \Phi)$. At the $E^2$-level, $E^2(\Phi)= \Delta_* \otimes
  {comp_{\omega}}_!$.}

Theorem 2 explain the interest to compute ${comp_{\omega}}_!$ in order
to do some computations of $\Phi$. So we have:

\noindent {\bf Theorem 3:} 

{\it The pointed loop-coproduct ${comp_{\omega}}_!$
  is zero.}

Theorem 2 and Theorem 3 proves that the loop-coproduct induced on the
Serre spectral sequence associated to the fibration $\Omega M \to LM
\stackrel{ev(0)}{\to} M$ vanishes at the $E^{\infty}$-level. In \cite{CT3}, the computation of this loop-coproduct is done for
$H_*(L \mathbb{C}P^n ; \mathbb{Q})$, using rationnal homotopy
theory. The two autors proves that in this case, the loop-coproduct is
non-zero. This indicates that the extension issues in our spectral
sequence are not trivial.

Now, let $N \to X \stackrel{p}{\to} M$ be a  locally 
trivial Serre fibration satisfying
hypothesis of proposition 1 of \cite{JF3} namely:

a) $N$, (respectively $M$)  is a finite dimensional smooth closed  oriented 
manifold of
 dimension $n$ (respectively $m$),

b) $M$ is a connected space and   $\pi_1(M)$ acts trivially on $H_\ast(N)$.

Then, we state the following theorem:

\noindent{\bf Theorem 4.}{\it Under the above hypothesis, the loop-coproduct
  $\Phi_X$ induces on the Serre spectral sequence associated to the fibration $LN \to LX  \stackrel{Lp}{\to} LM $ a structure of
coalgebra with a coproduct of degree $-(m+n) = dim X$. The tensor
product of coalgebra $(H_*(N), \Phi_N) \otimes (H_*(B), \Phi_B)$ is a
sub-coalgebra of $E^2[p]$.}

Let us recall the following theorem of \cite{JF3}:

\noindent {\bf Theorem B.} {\it Under hypothesis of 
 Proposition 1, the $(m,n)$-regraded Serre spectral sequence 
$
\{{\mathbb E}^r[Lp]\}_{r\geq 0}
$
 of the Serre fibration
$
LN \to LX \stackrel{Lp} \to LM
$
is a multiplicative spectral sequence which converges  to the algebra ${\mathbb 
H}_{\ast}(LX)$. Moreover the tensor product of graded algebras :
 $
{\mathbb H}_\ast (LM)\otimes {\mathbb H}_\ast (LN)$
 is a subalgebra of ${\mathbb E}^2[Lp]$.
In particular if $H_\ast (LM)$ is torsion free then
$
{\mathbb E}^2[p]= {\mathbb H}_\ast (LM)\otimes {\mathbb H}_\ast (LN)\,.
$}

As a immediate byproduct of theorem 4 and of this theorem B, we get:

\noindent{\bf Corollary.}{\it The regraded Serre spectral sequence associated
  to the fibration $
LN \to LX \stackrel{Lp} \to LM
$
is a spectral sequence of Frobenius algebra. Moreover the tensor
product of graded Frobenius algebras :
 $
{\mathbb H}_\ast (LM)\otimes {\mathbb H}_\ast (LN)$,
 is a sub-Frobenius algebra of ${\mathbb E}^2[Lp]$.}

\vspace{3mm}
{\bf Acknowledgement}
We would like to thank David Chataur and Jean-Claude Thomas for their
precious comments.

\vspace{3mm}

The paper is organized as follows:

{\bf 1) Definition of the loop coproduct

2) Proof of Proposition 1  

3)Proof of Theorem 1

4) Proof of Theorem 2

5) Proof of Theorem 3

6)  Proof of Theorem 4} 

\vspace{3mm}
\centerline{\bf \Large  1. Definition of the loop coproduct.}

Following \cite{CG3} and \cite{SUL3} the loop product (respectively the
loop coproduct) is a particular case of an "operation"
$$\mu_{\Sigma}: H_*(LM)^{\otimes p} \to H_*(LM)^{\otimes q} $$ 
with $p=2$, $q=1$ (respectively $p=1$, $q=2$) that will define a
positive boundary TQFT.
Here $\Sigma$ denotes an oriented surface of genus $0$ with a fixed
parametrization of the $p+q$ boundary components:
$$(S^1)^{\coprod p} \stackrel{in}{\to} \partial_{in} \Sigma
\hookrightarrow \Sigma \hookleftarrow \partial_{out} \Sigma
\stackrel{out}{\leftarrow} (S^1)^{\coprod q}$$
Applying the functor $Map(-,M)$ on gets the diagram:
$$(1) \qquad (LM)^{\times p} \stackrel{in^{\sharp}}{\leftarrow} Map(\Sigma,M)
\stackrel{out^{\sharp}}{\to} (LM)^{\times q}$$
Diagram (1) is homotopically equivalent to the diagram
$$(2) \qquad (LM)^{\times p} \stackrel{\rho_{in}}{\leftarrow} Map(c,M)
\stackrel{\rho_{out}}{\to} (LM)^{\times q}$$
where $c$ denotes a reduced Sullivan chord diagram with marking whitch
is associated to $\Sigma$ \cite{CG3}-§2.

In order to define the loop coproduct we restrict to the case $p=1$
and $q=2$. Then $\Sigma$ is the oriented surface of genus $0$ with one
incomming and two outcoming components of the boundary. The associated
Sullivan chord diagram with markings is determined by the pushout
diagram
$$
\xymatrix{
\ast \coprod \ast \ar[r] \ar[d] & \ast \ar[d] \\
S^1 \ar[r]^{in} & c }
$$
In particular, $c \simeq S^1 \vee S^1$, $S^1 \stackrel{in}{\to} c$ is
homotopic to the folding map $\nabla : S^1 \to S^1 \vee S^1$ and
$S^1 \coprod S^1 \stackrel{out}{\to}c$ is homotopic to the natural
projection $S^1 \coprod S^1 \to  S^1 \vee S^1$.
Diagram (2) then reduces to a commutative diagram of fibrations:
$$
\xymatrix{
LM \times LM \ar[d]^{ev(0) \times ev(0)} & LM \times_M LM
\ar[r]^{\delta_{in}} \ar[l]_{\delta_{out}} \ar[d]^{ev_{\infty}} & LM
\ar[d]^{ev(0) \times ev(1/2)}\\
M \times M & M \ar[l]^{\Delta}  \ar[r]^{\Delta} & M \times M }
$$ 
By definition (see \cite{CG3}), $\Phi = {\delta_{out}}_* \circ
{\delta_{in}}_!$.

\vspace{3mm}

\centerline{\bf \Large  2. Proof of Theorem 1}

Geometrically, $\delta_{out}$ is the inclusion of composable loops in
$LM \times LM$, namely $\tilde{\Delta}$. The other map, $\delta_{in}$,
can be thought as the inclusion of $map(8,M)$ in $LM$.  
Then this coproduct $\phi$ can be understood as follows:
the space $L_{1/2}M$ is the space of decomposable loops whitch embeds in $LM$.
This embedding is $d$-codimensionnal so that we can define the shriek
map of the embedding. Then, we decompose the loops of the space of decomposable loops.

The top line of the preceding diagram can be decomposed as follows:

$$ (*) \qquad LM \times LM \stackrel{\tilde{\Delta}=\delta_{out}}{\longleftarrow}
LM \times_M LM \stackrel{\gamma}{\to} L_{1/2}M
\stackrel{\tilde{j}}{\hookrightarrow} LM .$$
Observe that $\delta_{in}=\tilde{j} \circ \gamma = comp$ then $\Phi =
{\delta_{out}}_* \circ {\delta_{in}}_! = \tilde{\Delta}_* \circ
comp_!$. This proves
Theorem 1.

\rightline{$\square$}

\vspace{3mm}

\centerline{\bf \Large  3. Proof of Proposition 1}

By definition, $h_!$ is given by the following composition (see
\cite{JF3}):
$$ 
H_*(Y) \stackrel{inc_*}{\to} H_*(Y, Y - h(Y))
\stackrel{exc}{\to} H_*(\mbox{Tube} \, h / \partial \mbox{Tube} \, h)   \stackrel{exp_*}{\to} 
  H_*(\nu_Y, \partial \nu_Y)  
  \stackrel{\pi_*(\tau \cap -)}{\to} C_{*}(X)
$$
where $inc_*$ denotes the inclusion of pair, $exc$ is the excision
isomorphism, $exp$ is the homeomorphism given by the exponential
between the tubular neighbourood of the embedding and it's normal
bunddle and $\pi_*(\tau \cap - )$ , where $\pi$ denotes the projection
of the normal bunddle and $\tau \in H^*(\nu_Y, \partial \nu_Y)$ is the
Thom class of the embedding, is the Thom isomorphism. In the case we
consider, the normal bundle of $h$ is $\nu_Y =: * \to Y
\stackrel{h^{-1}}{\to} X $ . Then, $\mbox{Tube} \, h = Y = E(\nu_Y)$ and $\partial
\mbox{Tube} \, h = \emptyset = \partial E(\nu_Y)$ so that three first
application of the definition of $h_!$ are identity. We remark that
the projection of the normal bundle $\nu_Y$ is only $h^{-1}$ and that
$\tau$ lies in $H_0(\nu_Y, \partial \nu_Y) = H_0(Y)$. Then, the Thom
isomorphism is only $h^{-1}_*$, this achieve the proof Proposition 1.

\rightline{$\square$}

\vspace{3mm}

\centerline{\bf \Large  4. Proof of Theorem 2}
For the reader convenience, we recall here the definition of a fiber
embeding of \cite{JF3} and the main result of
\cite{JF3}.

\noindent{\bf Definition} We define a fiber embedding $(f,f^B)$ as a commutative diagram
$$ 
\begin{array}{cccc}
X&\stackrel{f} \to& X'\\
p\downarrow && \downarrow p'\\
B &\stackrel{f^B}\to& B'
\end{array}
$$
where
$$ (\ast)\quad \left\{
\begin{array}{lll}

a) & X, X',B \mbox{ and } B' \mbox{ are connected Hilbert  manifolds without 
boundary}\\

b) & f \mbox{ (respectively } f^B\mbox{) is a  smooth embedding of finite 
codimension } k_X \mbox{(respectively } 
k_B\mbox{)}\\

c)&p \mbox{ and } p' \mbox{ are locally trivial fibrations}\\

d)&\mbox{for some } b \in B \mbox{ the induced map}\\

&\hspace{32 mm} f^F : F:=p^{-1}(b)\to {p'}^{-1}(f(b)):=F'\\

&\mbox{ is an embedding of finite codimension } k_F \\

e)&\mbox{embeddings } f, f^B  \mbox{ and } f^F \mbox{admits Thom classes}
\,.
\end{array}
\right.
$$

\noindent{\bf First part of the main result.} {\it  Let $f: X \to X'$ be a fiber 
embedding as above.
 For each $n\geq 0 $ there exist filtrations
$\left\{ 
\begin{array}{lcl}

&\{0\}\subset F_0C_n(X)\subset F_1C_n(X)\subset ... \subset F_nC_n(X)=C_n(X)\\

&\{0\}\subset F_0C_n(X')\subset F_1C_n(X')\subset ... \subset F_nC_n(X')=C_n(X')
\end{array} 
\right.
$
and a chain representative  $ f_! : C_\ast(X')  {\to} C_{\ast-k_X} (X)\ $ of $
f_! : H_\ast(X')  {\to} H_{\ast-k_X} (X) $ satisfying:
$
f_!\left( F_\ast C_\ast (X')\right) \subset   F_{\ast-k_B} C_{\ast-k_X} (X')\,.
$}

Let $\{E^r[p]\}_{r\geq 0}  $ and $\{E^r[p']\}_{r\geq 0}$ be the
  spectral sequence induced by the above filtration.
\noindent{\bf Second part of the main result.} 
{\it The chain  
map $f_!$ induces a homomorphism of bidegree $(-k_B,-k_F)$ between the
associated  spectral sequences $\{ E^r(f_!)\} : \{E^r[p']\}_{r\geq 0}
\to \{E^r[p]\}_{r\geq 0}\,.$
There exists a chain representative   
$ f^B_! : C_\ast(B')  {\to} C_{\ast-{k_B}} (B)\mbox{ (respectively }
f^F_! : C_\ast(F') {\to} C_{\ast-k_F} (F)\mbox{)}$ of $ f^B_! :
H_\ast(B')  {\to} H_{\ast-k_B}(B) $ (respectively of $ f^F_! : 
H_\ast(F')  {\to}
H_{\ast-k_F} (F) $ ) such that
$$\{ E^2(f_!)\} = H_{\ast}(f^B_! ; {\cal H}_\ast(f^F_!)) : E^2_{s,t} [p'] = 
H_s(B' ;
{\cal H}_t(F')) \to  E^2_{s-{k_B}  ,t-{k_F}} [p] = H_{s-{k_B}}(B ; {\cal
H}_{t-{k_F}}(F))\,, $$
where  $ {\cal H}(-)$ denote the usual  system of local
coefficients. These spectral sequences are the Serre spectral
sequencesof the fibration.}

We use the description of the loop coproduct $\Phi$ given in the proof
of Theorem 1. Thus we obtain the following
commutative diagram of fibrations, where the central column is the
composite $(*)$:
$$
\xymatrix{
\Omega M \ar[r] & LM  \ar[r]^{ev(0)} & M \\
\Omega_{1/2} M  \ar[r] \ar[u]^{\tilde{i}}  \ar[d]^{\gamma_{\omega}^{-1}}& L_{1/2}M \ar[u]^{\tilde{j}}  \ar[d]^{\gamma^{-1}} \ar[r]^{ev(0)}  & M \ar[d]^= \ar[u]^=\\
\Omega M \times \Omega M  \ar[r]  \ar[d]^= & LM \times_M LM
\ar[d]^{\tilde{\Delta}}  \ar[r]^{ev_{\infty}} & M
\ar[d]^{\Delta} \\
\Omega M \times \Omega M  \ar[r]& LM \times LM  \ar[r]^{ev(0) \times ev(0)} & M \times M }
$$

We remark that $(\tilde{j} , id)$ is a fiber embedding in the sense of
\cite{JF3}, more precisely:

a)$LM$, $L_{1/2}M$, M are connected Hilbert manifolds without boundary.

b)$\tilde{j}$ (resp. $id$) is a smooth embedding of finite codimension d
(resp. $0$). 

c) $ev(0)$ and $ev(1/2)$ are locally trivial fibrations ($ev(1/2)$ is
locally trivial since it is homeomorphic to $ev_{\infty}$ whitch is
locally trivial by lemma 3 of \cite{JF3}).

d) $\tilde{i}$ is an embedding of finite codimension $d$.

e) $\tilde{j}$, $id$ and $\tilde{i}$ admit Thom classes. If we denote $\tau$ the
Thom class of $\tilde{j}$ and $u: \Omega M \to LM $ the canonical embedding,
then, the Thom class of $\tilde{i}$ is $u^*(\tau)$, \cite{CT3}.  

Then, we apply the main result of \cite{JF3} to prove that $\tilde{j}_!$
induces a morphism of spectral
sequences: $$E^*_{*,*}(\tilde{j}_!):E^*_{*,*}[ev(0)] \to
:E^*_{*-d,*}[ev(1/2)].$$
By naturality of the \sss s, $\gamma^{-1}$ and $\tilde{\Delta}$ induce a
morphism  of spectral sequences:$$E^*_{*,*}(\gamma^{-1}_*):E^*_{*-d,*}[ev(1/2)] \to
E^*_{*-d,*}[ev_{\infty}]$$ and :$$E^*_{*,*}(\tilde{\Delta}_*):E^*_{*-d,*}[ev_{\infty}] \to
E^*_{*-d,*}[ev(0) \times ev(0)].$$
Then, composing this morphisms, we define a coproduct $$E^*(\Phi):
E^*_{*,*}[ev(0)] \to E^*_{*-d,*}[ev(0) \times ev(0)]$$ induced by
$\Phi$. On the base of the fibration $ev(0)$, $\Phi$ induces
$\Delta_*$. On the fiber, $\Phi$ induces ${comp_{\omega}}_!$ such
that at the $E^2$-level, $E^2(\Phi) = \Delta_* \otimes
{comp_{\omega}}_!$.

\rightline{$\square$}

\vspace{3mm}

\centerline{\bf \Large  5. Proof of Theorem 3}

Consider the following pull-back diagram:
$$
\xymatrix{
\Omega_{1/2} M \ar[r]^= \ar[d] & \Omega_{1/2} M  \ar[d] \\  
\Omega_{1/2} M \ar[r]^{\tilde{i}} \ar[d] & \Omega M \ar[d]^{ev(1/2)} \\
\ast \ar[r]^i & M }
$$
where $i$ and $\tilde{i}$ are defined in the introduction.
From Lemma 3 of \cite{JF3}, this diagram is a fiber embedding thus
$\tilde{i}_!$ induces a morphism of \sss s : $E^*(\tilde{i}_!):E^*_{*,*}(2) \to
E^*_{*-d,*}(1)$ where $E^*_{*,*}(1)$ (resp. $E^*_{*,*}(2)$) denotes
the \sss $\,$ 
associated to the left fibration of the above diagram (resp.
the \sss $\,$ 
associated to the right fibration of the above diagram). 
Moreover, $\tilde{i}_*$ induces a morphism of \sss s \cite{MC3}:
 $E^*(\tilde{i}_*):E^*_{*,*}(1) \to E^*_{*,*}(2)$.
We know that $\tilde{i}_* \circ \gamma_* = comp_*$ is onto ($comp_*$
has a unit). Since $\gamma_*$ is an isomorphism, this proves that
$\tilde{i}$ is onto. The spectral sequence  $E^*_{*,*}(1)$ collapses at the
$E^2$-level, it has only one column at this level: $$E^n_{0,*}(1)
\simeq H_*(\Omega_{1/2} M) \, , \, n \geq 2.$$
Then, since $\tilde{i}_*$ is onto, $E^{\infty}(\tilde{i}_*)$ is onto.
Moreover, $E^2(\tilde{i}_*) = i_* \otimes id$ thus $Im(E^2(\tilde{i}_*)) =
E^2_{0,*}(2)$ so that $Im(E^2(\tilde{i}_*)) \subset  E^{\infty}_{0,*}(2)$
that is why $ E^{\infty}_{*,*}(2) = Im(E^{\infty}(\tilde{i}_*)) \subset
E^{\infty}_{0,*}(2)$. 

We have proved that $$E^{\infty}_{*,*}(2)=E^{\infty}_{0,*}(2).$$
But  $E^2(\tilde{i}_!) = i_! \otimes id$ and $\tilde{i}_!$ is non zero only in
degree $d$ with values in degree $0$ ($i_!$ send the fundemental class
of $H_d(M)$ on a generator of $H_0(M)$ and is zero elsewhere). Then,
at the $E^2$-level of $E^*_{*,*}(2)$, only the column of abcisse $d$
has a non zero image by $E^2(\tilde{i}_!)$. Consequently, at the
aboutment, only  $E^{\infty}_{d,*}(2)$ has a non zero image by
$E^{\infty}(\tilde{i}_!)$. We have shown that  $E^{\infty}_{d,*}(2)=0$
thus $Im(E^{\infty}(\tilde{i}_!))=0$.

\rightline{$\square$} 

\vspace{3mm}

\centerline{\bf \Large  6. Proof of Theorem 4}
{\bf Notations.} In this section, we use the same notations as in the
introduction but we add a subscript to indicate which manifold we
refer (for example, $\Phi_X$denotes the loop-coproduct on $H_*(LX)$).

First, we need the following lemma:

{\bf Lemma.}{\it The commutative diagram 
$\xymatrix{
LN \ar[r] & LX \ar[r] & LM \\
L_{1/2}N \ar[r] \ar[u]& L_{1/2}X \ar[r] \ar[u] & L_{1/2}M \ar[u] }$
is a fiber embedding}

Proof of the Lemma. We check that the diagram verify all the
properties of the definition of a fiber embedding. The point a) of the definition is verified in proof of Theorem 2. The same holds
for the point b) with $i_X$ and $i_M$ smooth embedding of codimension
$x = m+n$ and $m$. The point c) comes from the fact that the fibration  
$\xymatrix{
LN \times_N LN \ar[r] & LX \times_X LX \ar[r]^{Lp} & LM \times_M LM }$
is locally trivial (see \cite{JF3} 3.3) and since it is homeomorphic
to the fibration  $\xymatrix{
L_{1/2}N \ar[r] & L_{1/2}X  \ar[r] & L_{1/2}M  }$, this last fibration
is locally trivial.
The point d) comes directly from the diagram of the Lemma. Since their
exists Thom class for the embedding $L_{1/2}M
\hookrightarrow LM$ (see proof of theorem 2 part 4), the embeddings  $L_{1/2}N
\hookrightarrow LN$,  $L_{1/2}X
\hookrightarrow LX$ and  $L_{1/2}M
\hookrightarrow LM$ admit Thom class. This proves the point e).

\rightline{$\square$}

Proof of Theorem 4:
we consider the following commutative diagram:
$$
\xymatrix{
LN \ar[r] & LX \ar[r] & LM \\
L_{1/2}N \ar[r] \ar[u]^{\tilde{j}_N} & L_{1/2}X  \ar[r] \ar[u]^{\tilde{j}_X}& L_{1/2}M
\ar[u]^{\tilde{j}_M}\\
LN \times_N LN \ar[r] \ar[u]^{h_N} \ar[d]_{\tilde{\Delta}_N} & LX
\times_X LX \ar[r] \ar[u]^{h_X} \ar[d]_{\tilde{\Delta}_X} & LM \times_M LM  \ar[u]^{h_X}
\ar[d]_{\tilde{\Delta}_M} \\
LN \times LN \ar[r] & LX \times_X LX \ar[r] & LM \times_M LM }.$$

It comes from the Lemma and the main result of \cite{JF3} that
$\tilde{j}_{X!}$ induces a morphism of spectral sequence of degree $-x$ so
that $\Phi_X= \tilde{\Delta}_{X*} \circ h_{X*} \circ i_{X!}$ induces a
morphism of spectral sequence $E^*[\Phi_X]$. This morphism induces on
the homology of $M$  $\Phi_M = \tilde{\Delta}_{M*} \circ h_{M*} \circ
\tilde{j}_{M!}$ and on the homology of $N$  $\Phi_N = \tilde{\Delta}_{N*}
\circ h_{N*} \circ \tilde{j}_{N!}$. This achieve the proof of Theorem 4.

\rightline{$\square$}

\end{document}